\documentstyle[12pt,amssymb]{article}

\begin{document}

\newtheorem{thm}{Theorem}[section]
\newtheorem{lemma}[thm]{Lemma}
\newtheorem{defn}{Definition}[section]
\newtheorem{prop}[thm]{Proposition}
\newtheorem{corollary}[thm]{Corollary}
\newtheorem{remark}{Remark}

\renewcommand{\theequation}{\arabic{section}.\arabic{equation}}

\catcode`\@=11
\renewcommand{\section}{
        \setcounter{equation}{0}
        \@startsection {section}{1}{\z@}{-3.5ex plus -1ex minus
        -.2ex}{2.3ex plus .2ex}{\Large\bf}
}
\catcode`\@=12

\def\ee{\varepsilon}
\def\qed{\hfill\hbox{\vrule height1.5ex width.5em}}
\def\RR{{\bf R}}
\def\KK{{\bf K}}
\def\pf{\noindent{\bf Proof.} } 

\title{\Large \bf Martin Boundary and Integral Representation for 
Harmonic Functions of Symmetric Stable Processes}
\author{Zhen-Qing Chen
\thanks{The research of this author is supported in part by an NSA grant}\\
Departments of Mathematics and Statistics\\ 
Cornell University \\
 Ithaca, NY 14853, USA \\
Email: zchen@math.cornell.edu \smallskip \\
and  
\smallskip \\
Renming Song\\
Department of Mathematics\\
University of Illinois \\
Urbana, IL 61801, USA \\
Email: rsong@math.uiuc.edu }
\date{ }
\maketitle

\begin{abstract}
Martin boundaries and integral representations of positive
functions which are harmonic in a bounded domain $D$ with respect to
Brownian motion are well understood. Unlike the Brownian case,
there are two different kinds of harmonicity with respect to a 
discontinuous symmetric
stable process. One kind are functions harmonic in $D$ with respect
to the whole process $X$, and the other are functions harmonic
in $D$ with respect to the process $X^D$ killed upon leaving
$D$. In this paper we show that for bounded Lipschitz domains, the
Martin boundary with respect to the killed stable process
$X^D$ can be identified with
the Euclidean boundary.  We further give integral representations
for both kinds of positive harmonic functions. Also given 
is the conditional
gauge theorem conditioned according to Martin kernels and
the limiting behaviors of the $h$-conditional stable process,
where $h$ is a positive harmonic function of $X^D$.
In the case when
$D$ is a bounded $C^{1, 1}$ domain, sharp estimate on the Martin
kernel of $D$ is obtained.
\end{abstract}

\vspace{.3in}


\vspace{.1in}
\noindent {\bf Running Title:} Martin Boundary for Stable Processes

\pagebreak

\begin{doublespace}

\section{Introduction}

Martin boundary and integral representation for harmonic functions
of diffusions processes (or of elliptic differential operators) 
are well studied. However there is little detailed analysis
of these for Markov processes with jumps (or for integro-differential
operators). 
In this paper we take a closer look at an important class
of discontinuous Markov processes---symmetric $\alpha$-stable processes
with $0< \alpha <2$, and study the notion and integral representation
of harmonic functions for these processes, where some new phenomena arise.
We hope that this paper can shed some new light on the potential theory
of general Markov processes.

Symmetric stable processes constitute an important
 subfamily  of L\'evy processes.
A symmetric $\alpha$-stable process $X$ on $\RR^n$
is a L\'evy process
whose transition density $p(t, x-y)$
relative to the Lebesgue measure is uniquely determined
by its Fourier transform $\int_{\RR^n} e^{i x\cdot \xi}  p(t, x) dx
=e^{-t |\xi|^\alpha} $. Here $\alpha$ must be in the interval $(0, \, 2]$.
When $\alpha=2$, we get a Brownian motion running with a time clock
twice as fast as  the standard one.
In this paper, 
symmetric stable processes are referred to the case
when $0<\alpha<2$, unless otherwise specified.

Unlike the Brownian case,
there are two different kinds of harmonicity with respect to symmetric
stable processes, one kind are functions harmonic in $D$ with respect
to the whole process $X$, and the other are functions harmonic
in $D$ with respect to the process $X^D$ killed upon leaving
$D$. The theory of Martin kernel and Martin boundary for
the killed process $X^D$ is known from the general theory.
This Martin boundary gives an integral representation
for positive functions harmonic in a domain $D$ with respect
to the killed process $X^D$. 
We show that when $D$ is a bounded Lipschitz domain,
 the Martin boundary with respect to killed symmetric stable process 
 $X^D$ in $D$
 coincide with the Euclidean boundary.
It seems that integral representations
of positive  functions harmonic in a domain $D$ with respect
to the whole processes $X$  have not been studied in the literature.
In this paper, we present an integral
representation for positive functions harmonic in a domain $D$ with respect
to the whole processes and this representation is shown to be unique.
In particular, this implies that any harmonic function
with respect to the whole process is uniquely determined 
by its values in $D$.
In the case when $D$ is a bounded $C^{1, 1}$ domain, 
sharp estimates on the Martin kernel are given. As a consequence of
these estimates, we prove a conditional
gauge theorem conditioned according to  Martin kernel.
We also study
the limiting behavior of the $h$--conditioned symmetric stable process
in $D$ when $h$ is a positive harmonic function of $X^D$, and the 
the limiting behavior of the $h$--conditioned symmetric stable process
will provide a probabilistic interpretation to positive
harmonic functions of $X^D$.

This paper is organized as follows. The definitions
of harmonic and superharmonic functions with respect to symmetric
stable processes are given in section 2. 
Some important facts about those harmonic functions are also given in
section 2.
Section 3 contains results on Martin boundary and conditional gauge theorem.
Integral representations
of positive  functions harmonic in a domain $D$ with respect
to the whole processes are given in section 4.

\medskip

In the sequel, we will use $v^+$ and $v^-$ to
denote the positive and negative part of
a real-valued Borel measurable function $v$, i.e.,  $v^+=\max\{ v ,\, 0\}$
and $v^- =\max \{-v ,\, 0\}$.

\vspace{.1in}

\noindent{\bf Acknowledgement.}
After the first draft of this paper was written, we met Krzysztof Bogdan
and learned about his paper \cite{bog2} at a 
Stochastic Analysis workshop held at the MSRI at Berkeley
from March 22--27, 1998.
Some of our results had been independently and simultaneously
obtained by him in \cite{bog2}. However
the approach of this paper and that of 
\cite{bog2} are different.
We  thank Krzysztof Bogdan
for the very interesting and helpful discussions
at the MSRI, Berkeley and at University of Washington, Seattle.
We are also grateful to Chris Burdzy, Eugene Dynkin, Pat Fitzsimmons,
Tom Kurtz and Sergei Kuznetsov for very helpful discussions about this paper.
Parts of the research for this paper were conducted while the authors
were visiting the CMS at University of Wisconsin at Madison and the MSRI in
Berkeley. Financial support from  two institutions are gratefully
acknowledged.

\section{Definitions and Preliminaries}

In sections 2--4 of this paper, we always assume $n\geq 2$.
Let $X=(\Omega, {\cal F}, X_t, {\cal F}_t, P_x)$ 
be a symmetric $\alpha$-stable process on $\RR^n$ with
$0 < \alpha < 2$, where $\{ {\cal F}_t, \, t\geq 0\}$
is the minimal admissible $\sigma$-fields generated by $X$. 
The process $X$ is transient and we are going to use
$G$ to denote the potential of $X$.
We know that the Green function of $X$ is given by
$$
  G(x, y)=
2^{-\alpha} \pi^{-{n\over2 }}\Gamma \left({n-\alpha\over 2}\right)
\Gamma\left({\alpha \over 2}\right)^{-1}  
 |x-y|^{\alpha -n}.   
$$
For a domain $D$ in $\RR^n$, let $\tau_D =\inf\{ t>0: X_t \not\in D \}$.
Adjoin a cemetery point $\partial$ to $D$ and set
$$
X^D_t(\omega)
   = \cases  {X_t(\omega) \ & if $t< \tau_D$, \cr
     \partial \   & if $t \ge \tau_D$.}
$$
$X^D$ is a strong Markov process with state space $D_{\partial}=
D\cup\{\partial\}$,
which is called the subprocess of the 
symmetric $\alpha$--stable process $X$  killed upon leaving $D$,
or simply the symmetric $\alpha$--stable process in $D$.
We are going to use $G_D$ to denote the Green function
of $X^D$.

For Brownian motion or other
diffusion processes, there is only one kind of harmonicity on a
domain $D$. However, for symmetric stable processes, there are two
kinds of harmonic functions on $D$: functions which are harmonic in
$D$ with respect to the killed process $X^D$ and functions which
are harmonic in $D$ with respect to the process $X$. The precise
definitions of these two kinds of harmonic functions are as follows.

\begin{defn}\label{def:har1}
Let $D$ be a domain in $\RR^n$. 
A locally integrable function
$f$ defined on $D$ taking values in $(-\infty, \, \infty]$ and
satisfying the condition $\int_{\{|x|>1\}\cap D} |f(x)| |x|^{-(n+\alpha )}
dx <\infty$
is said to be 

\begin{description}
\item 1) harmonic with respect to $X^D$ if 
$f$ is continuous in $D$ and for each $x \in D$ and
each ball $B(x, r)$ with $\overline{B(x, r)}\subset D$,
$$
f(x)= E_x[f(X_{\tau_{B(x, r)}}); \tau_{B(x, r)}<\tau_D];
$$

\item 2) superharmonic respect to $X^D$ if 
$f$ is lower semicontinuous in $D$ and for each $x \in D$ and
each ball $B(x, r)$ with $\overline{B(x, r)}\subset D$,
$$
f(x)\ge E_x[f(X_{\tau_{B(x, r)}}); \tau_{B(x, r)}<\tau_D].
$$
\end{description}
\end{defn}

The next definition is taken from Landkof \cite{ns:lan}.

\begin{defn}\label{def:har2}
Let $D$ be a domain in $\RR^n$. 
A locally integrable function
$f$ defined on $\RR^n$ taking values in $(-\infty, \, \infty]$ and
satisfying the condition $\int_{\{|x|>1\}} |f(x)| |x|^{-(n+\alpha )}
dx <\infty$
is said to be 

\begin{description}
\item 1) harmonic in $D$ with respect to $X$ if 
$f$ is continuous in $D$ and for each $x \in D$ and
each ball $B(x, r)$ with $\overline{B(x, r)}\subset D$,
$$
f(x)= E_x[f(X_{\tau_{B(x, r)}})];
$$

\item 2) superharmonic in $D$ with respect to $X$ if 
$f$ is lower semicontinuous in $D$ and for each $x \in D$ and
each ball $B(x, r)$ with $\overline{B(x, r)}\subset D$,
$$
f(x)\ge E_x[f(X_{\tau_{B(x, r)}})].
$$
\end{description}
\end{defn}

\noindent{\bf Remark 2.1} (1) If $f$ is a lower semicontinuous
function defined on $D$ taking values in $(-\infty , \, \infty]$,
then $f$ is bounded from below on any subdomain whose closure 
is contained in $D$. Thus for such kind of function $f$
which is locally integrable and satisfying \break
$\int_{\{|x|>1\}} |f(x)| |x|^{-(n+\alpha )}
dx <\infty$, it follows from estimate (\ref{eqn:2.2}) below that
$ E_x[f^-(X_{\tau_{B(x, r)}})] <\infty$ for any ball
$B(x, r)$ with $\overline{B(x, r)}\subset D$. Therefore the expections
in Definitions \ref{def:har1} and \ref{def:har2} are well defined.

\smallskip

(2) For a function $f$ which is (super)harmonic with respect to $X^D$, if
we extend it to be zero off the domain $D$, then the resulting function 
is (super)harmonic in $D$ with respect to $X$.

\smallskip

(3) Conversely, if $f$ is a non-negative superharmonic in $D$ 
with respect to $X$, then clearly it is a superharmonic with respect to $X^D$.

\medskip

We now record some facts, which will be used later,
concerning bounded Lipschitz
domains and  the exit distributions
of $X$ from a domain $U$.
Recall that 
a bounded domain $D$ in $\RR^n$ is said to be 
a bounded Lipschitz domain with
Lipschitz characteristic constants
$(r_0, \, A_0)$ if for every $z \in \partial D$, there is a local 
coordinate system $(\xi_1,\xi^{(1)}) \in \RR \times \RR^{n-1}$
 with origin sitting at $z$ and there is a Lipschitz function $f$
defined on $\RR^{n-1}$ with
Lipschitz constant $A_0$ such that
$D \cap B(z,r_0) = B(z,r_0) \cap \{\xi = (\xi_1,\xi^{(1)}):
\xi_1 > f(\xi^{(1)})\}$.
A domain $U$ in $\RR^n$ is said to satisfy the {\it uniform
exterior cone condition} if there exist constants $\eta>0$,
$r>0$ and a cone
${\cal C} =\{x=(x_1, \dots, x_n)\in \RR^n: 0<x_n,
(x^2_1+\cdots+x^2_{n-1})^
{1/2}<\eta x_n\}$ such that for every $z\in \partial U$, there is
a cone ${\cal C}_z$ with vertex $z$, isometric to ${\cal C}$ and
satisfying ${\cal C}_z\cap B(z, r)\subset U^c$.
It is well known that bounded Lipschitz domains satisfy 
the uniform exterior cone condition. 
Bogdan showed in \cite{kr:bog} that for a 
bounded domain $U$ satisfying  the uniform
exterior cone condition, 
\begin{equation}\label{eqn:2.1}
P_x (X_{\tau_U} \in \partial U) =0 \quad \mbox{for all } x \in U.
\end{equation}

In \cite{s1:che}, Chen and Song showed that if
$U$ is a bounded $C^{1,1}$ domain in $\RR^n$, there is a 
Poisson kernel $K_U(x, z)$ defined on $U\times (\RR^n \setminus
\overline U)$ such that for any bounded Borel measurable function
$\phi$, $E_x[ \phi(X_{\tau_U})] =\int_{U^c} \phi(z) K_U(x, z) dz$,
where $dz$ is the Lebesgue measure on $\RR^n$. Furthermore
there  exists a $C=C( U, \alpha)>1$ such that
for $ x\in U$ and $ z\in {\overline U}^c$,
\begin{equation}\label{eqn:2.2}
 \frac{\delta(x)^{\alpha/2}}
        {C\, \delta(z)^{\alpha/2}(1+\delta(z))^{\alpha/2}} \
       \frac{1}{|x-z|^n}
   \le K_U(x,z)
   \le  \frac{C \, \delta(x)^{\alpha/2}}
        {\delta(z)^{\alpha/2}(1+\delta(z))^{\alpha/2}} \
       \frac{1}{|x-z|^n},
\end{equation}
where $\delta (y)=\hbox{dist}(y, \partial U)$ is the Euclidean 
distance from point $y$ to the set $\partial U$.
Here a domain $U$ is $C^{1, 1}$ means that for every $z\in \partial D$,
there exists a $r>0$ such that $B(z, r)\cap\partial D$
is the graph of a function whose first derivatives are Lipschitz.

\medskip

It is well known that for any domain $D$, there exists 
an increasing sequence of bounded $C^\infty$-smooth domains
$\{D_k \}_{k \geq 1}$ such that $\overline {D_k} \subset D_{k +1}$
for $k\geq 1$ and that $\bigcup_{k=1}^\infty D_k =D$ (see, for example,
Lemma 2.4 of \cite{chen}).

\medskip

\begin{thm}\label{1a}
Suppose that $D$ is a bounded domain in $\RR^n$.
If $h$ is  superharmonic in $D$ with respect to $X$, then for
any domain $D_1\subset\overline{D_1}\subset D$,
$E_x[ h^-(X_{\tau_{D_1}})] <\infty$ and
\[
h(x)\ge E_x[h(X_{\tau_{D_1}})] \quad \mbox{for every } x\in D_1 .
\]
\end{thm}

\pf 
For a fixed $\epsilon>0$ and each $x\in D_1$ we put
\[
r(x)=\frac12\delta(x, \partial D_1)\wedge\epsilon,
\,\,\, B(x)=B(x, r(x)),
\]
where $\delta(x, \partial D_1)$ denotes the Euclidean distance
between $x$ and $\partial D_1$.
 Define a sequence of stopping times $\{T_m,
m\ge1\}$ as follows:
\[
T_1=\inf\{t>0: X(t)\notin B(X_0)\},
\]
and for $m\ge 2$,
\[
T_m=\cases{ T_{m-1}+
T_{B(X_{T_{m-1}})}\circ\theta_{T_{m-1}} & if $X_{T_{m-1}})
\in D_1$, \cr
\tau_{D_1} & otherwise. \cr}
\]
The superharmonicity of $h$ and the strong Markov property imply that
\[
h(X_{T_{m-1}})\ge E_x[h(X_{T_m})|{\cal F}_{T_{m-1}}].
\]
Thus $\{h(X_{T_m}), m\ge1\}$ is a supermartingale under $P_x$.

We claim that for each $x\in D_1$,
\[
P_x(\lim_{m\to \infty}T_m=\tau_{D_1})=1.
\]
It is clear that $P_x$--a.\/s., $T_m\uparrow$ and 
$T_m\le\tau_{D_1}<\infty$ because $D$ is bounded. Let
$T_{\infty}=\lim_{m\to\infty}T_m$. Then $X_{T_{\infty}}=
\lim_{m\to\infty}X_{T_m}$ by quasi left continuity. On
the set $\{T_{\infty}<\tau_{D_1}\}$,
we have $X_{T_{\infty}}\in D_1$. Then for all sufficiently
large values of $m$, we have 
$\delta(X_{T_{m-1}}, \partial D_1)>\frac12\delta(X_{T_{\infty}},
\partial D_1)>0$ and $|X_{T_{m-1}}-X_{T_m}|<\frac14
\delta(X_{T_{\infty}},\partial D_1)\wedge\epsilon$. But by
the definition of $T_m$, $|X_{T_{m-1}}-X_{T_m}|
>\frac12\delta(X_{T_{m-1}}, \partial D_1)\wedge\epsilon$. 
These inequalities are incompatible. Hence
$P_x(T_{\infty}<\tau_{D_1})=0$.

Put $A=\{ \tau_{D_1}=T_m \mbox{ for some } t\ge1\}$.
Since $h$ is bounded from below on $\overline{D_1}$, we have by
Fatou's lemma
\[
\liminf_{m\to\infty}E_x[h(X_{T_m}); T_m<\tau_{D_1}]
\ge E_x[\liminf_{m\to\infty}h(X_{T_m}); A^c]\\
\ge E_x[h(X_{\tau_{D_1}}); A^c].
\]

Take two smooth domains $D_2$ and $D_3$ such that 
$\overline{D_1}\subset D_2\subset\overline{D_2}\subset D_3\subset\overline{D_3}
\subset D$, then $h$ is bounded from below
on $\overline{D_3}$.
Since
\begin{eqnarray*}
h(X_{\tau_{D_1} })&=& h(X_{\tau_D})1_{\{ X_{\tau_{D_1}}\in
D_3 \} } + h(X_{\tau_{D_1} }) 1_{ \{X_{\tau_{D_1} }\notin
 D_3 \} }\\
&= & h (X_{\tau_{D_1}}) 1_{\{ X_{\tau_{D_1} }\in
 D_3 \} } + h (X_{\tau_{D_2}}) 1_{\{ X_{\tau_{D_1}}\notin
  D_3 \}},
\end{eqnarray*}
we have by estimate (\ref{eqn:2.2}) with $D_2$ in place of $U$ 
and the integrability assumption about $h$ in
Definition \ref{def:har2} that
\[
E_x[ h^- (X_{\tau_{D_1} })]<\infty.
\]
Thus by Fatou's Lemma  
\begin{eqnarray*}
h(x)&\geq &\liminf_{m\to\infty}E_x[h(X_{T_m})]\\
&\ge&\liminf_{m\to\infty}E_x[h(X_{\tau_{D_1}})
, T_m=\tau_{D_1}]
+\liminf_{m\to\infty}E_x[h(X_{T_m }), T_m<\tau_{D_1}]\\
&\ge&E_x[h(X_{\tau_{D_1}}); A]+E_x[h(X_{\tau_{D_1}}); A^c]\\
&= & E_x[h(X_{\tau_{D_1}})].
\end{eqnarray*}
This completes the proof.
\qed

\vspace{.1in}

\begin{thm}\label{1b}
Suppose that $D$ is a bounded domain in $\RR^n$.
If $h$ is a harmonic in $D$ with respect to $X$, then for
any domain $D_1\subset\overline{D_1}\subset D$,
$h(X_{\tau_{D_1}})$ is $P_x$-integrable and 
\begin{equation}\label{eqn:har}
h(x)= E_x[h(X_{\tau_{D_1}})] \quad \mbox{for every }x\in D_1 .
\end{equation}
\end{thm}

\pf
We can always take a smooth domain $D_2$ such that
$D_1\subset\overline{D_1}\subset D_2\subset\overline{D_2}\subset D$.
If we could prove that for any $x\in D_2$, $h(X_{\tau_{D_2}})$
is $P_x$-integrable and
\[
h(x)= E_x[h(X_{\tau_{D_2}})].
\]
then by strong Markov property we immediately get
$h(X_{\tau_{D_1}})$ is $P_x$-integrable and 
\[
h(x)= E_x[h(X_{\tau_{D_1}})],\qquad x\in D_1.
\]
Therefore we can assume, without loss of generality, that
$D_1$ is a smooth domain.

Define $T_m$ as in the proof of the previous theorem, then in this case
$\{h(X_{T_m}), m\ge1\}$ is a martingale under $P_x$ for any $x\in D_1$.
By (\ref{eqn:2.1}) with $D_1$ in place of $U$, we have
 $P_x(\tau_{D_1}=T_m\, \mbox{ for some $m\ge1$})=1$.
Since $h$ is bounded on $D_1$, we have
\[
|E_x[h(X_{T_m}), T_m<\tau_{D_1}]|\le C
P_x(T_m<\tau_{D_1})\rightarrow 0.
\]
Take a domain $D_2$ such that $\overline{D_1}\subset D_2
\subset\overline{D_2}\subset D$, then $h$ is continuous
and therefore bounded on $\overline{D_2}$.
By the estimate (\ref{eqn:2.2}) with $D_1$ in place of $U$ and 
the integrability assumption about $h$ in Definition 2.2 
we have
$E_x[|h|(X_{\tau_{D_1}})]
<\infty$. Thus by the dominated convergence theorem
\[
\lim_{m\to\infty}E_x[h(X_{\tau_{D_1}})
, T_m=\tau_{D_1}]=E_x[h(X_{\tau_{D_1}})].
\]
Therefore
\begin{eqnarray*}
h(x)& =& \lim_{m\to\infty}E_x[h(X_{T_m})] \\
&=& \lim_{m\to\infty}E_x[h(X_{\tau_{D_1}})
, T_m=\tau_{D_1}]
+ \lim_{m\to\infty} E_x[h(X_{T_m}), T_m<\tau_{D_1}]\\
& = & E_x[h(X_{\tau_{D_1}})].
\end{eqnarray*}
\qed

\vspace{.1in}

Similarly, we have the following result for functions
harmonic with respect to $X^D$.

\begin{thm}\label{1c}
Suppose that $D$ is a bounded domain in $\RR^n$.
If $h$ is superharmonic in $D$ with respect to $X^D$, then for
any domain $D_1\subset\overline{D_1}\subset D$,
$E_x [h^- (X^D_{\tau_{D_1}})] <\infty$ and 
\[
h(x)\ge E_x[h(X^D_{\tau_{D_1}})] \quad \mbox{for every } x\in D_1 .
\]
If $h$ is harmonic in $D$ with respect to $X^D$, then for
any domain $D_1\subset\overline{D_1}\subset D$,
$h(X^D_{\tau_{D_1}})$ is $P_x$-integrable and 
\[
h(x)= E_x[h(X^D_{\tau_{D_1}})] \quad \mbox{for every } x\in D_1.
\]
\end{thm}

\begin{thm}\label{1d}
Suppose that $D$ is a bounded  domain and $h$ is harmonic
in $D$ with respect to $X$ and continuous on
$\overline{D}$, then $h(X_{\tau_D})$ is $P_x$-integrable and 
\[
h(x)=E_x[h(X_{\tau_D})], \qquad \hbox{for each } x\in D.
\]
\end{thm}

\pf
Take an increasing sequence of smooth
domains $\{D_m\}_{m\geq 1}$ such that $
\overline{D_m}\subset D_{m+1}$ and
$\bigcup_{m=1}^\infty D_m =D$.
Set $\tau_m=\tau_{D_m}$. Then $\tau_m\uparrow
\tau_D$ and $\lim_{m\to\infty}X_{\tau_m}=X_{\tau_D}$ by quasi-left 
continuity of $X$. Set
\begin{equation}\label{eqn:2.3}
A=\{\tau_m=\tau_D \mbox{ for some } m\ge1\}.
\end{equation}
From Theorem \ref{1b} we know that for any $m\ge1$,
\[
h(x)=E_x[h(X_{\tau_m})], \qquad x\in D_m.
\]
Since $h$ is continuous on $\overline{D}$, we have 
by dominated convergence theorem that
\begin{equation}\label{eqn:2.4}
\lim_{m\to\infty}E_x[h(X_{\tau_m}), \tau_m<\tau_D]
=E_x[h(X_{\tau_D}), A^c].
\end{equation}
Since $h$ is continuous on $\overline{D}$, we can find
two smooth domains $U_1$ and $U_2$ such that
$\overline{D}\subset U_1\subset
\overline{U_1}\subset U_2$ and that $h$ is bounded
on $\overline{U_2}$.
Since
\begin{eqnarray*}
|h|(X_{\tau_D})&=& |h|(X_{\tau_D})1_{\{ X_{\tau_D}\in
U_2 \} } + |h|(X_{\tau_D}) 1_{ \{X_{\tau_D}\notin
 U_2 \} }\\
& =  & |h| (X_{\tau_D}) 1_{\{ X_{\tau_D}\in
 U_2\} } + |h| (X_{\tau_{U_1}}) 1_{\{ X_{\tau_{U_1}}\notin
  U_2 \}},
\end{eqnarray*}
we have by estimate (\ref{eqn:2.2}) with $U_1$ in place of $U$ that
\[
E_x[|h|(X_{\tau_D})]<\infty.
\]
Thus by the dominated convergence theorem
\begin{eqnarray*}
h(x)&=&\lim_{m\to\infty}E_x[h(X_{\tau_m})]\\
&=&\lim_{m\to\infty}E_x[h(X_{\tau_D})
, \tau_m=\tau_{D}]
+\lim_{m\to\infty}E_x[h(X_{\tau_m}), \tau_m<\tau_D]\\
&=& E_x[h(X_{\tau_D})].
\end{eqnarray*}
\qed

\medskip

\noindent{\bf Remark 2.2} If $D$ is a bounded domain satisfying
the uniform exterior cone condition, then the conclusion of
Theorem \ref{1d} holds
for any harmonic function $h$ in $D$ with respect to $X$ that is
bounded in a neighborhood of $\overline D$. This is because in this
case by (\ref{eqn:2.1}) $P_x (A) =1$ for $x\in D$, where $A$ is the
set defined in (\ref{eqn:2.3}) and the term in (\ref{eqn:2.4}) vanishes.
The rest of the argument goes through without the continuous assumption
on $h$ up to the boundary $\partial D$. 

\vspace{.1in}

Obviously there are plenty of bounded functions which
are harmonic in $D$ with respect to the whole processes $X$. 
The following results says that, when $D$ is a bounded
domain satisfying the uniform exterior cone condition,
the only bounded 
function which is harmonic in $D$ 
with respect to $X^D$ is constant zero.

\begin{thm}\label{1e}
Suppose that $D$ is a bounded  domain in $\RR^n$
satisfying the uniform exterior cone condition. If
$h$ is a bounded function harmonic in $D$ with respect to $X^D$, 
then $h$ must be identically zero.
\end{thm}

\pf
Take an increasing sequence of smooth domains $D_m$ such that $D_m\subset
\overline{D_m}\subset D_{m+1}\subset\overline{D_{m+1}}
\subset D$ and set $\tau_m=\tau_{D_m}$. Then $\tau_m\uparrow
\tau_D$. By (\ref{eqn:2.1}), 
we know that $P_x(\tau_D=\tau_m\, \mbox{ for some $m\ge 1$})
=1$ for $x\in D$. From Theorem \ref{1b} we get that
\begin{eqnarray*}
|h(x)|&=&|E_x[h(X^D_{\tau_m})]|\\
&\le& C\, P_x(\tau_m<\tau_D)\rightarrow0.
\end{eqnarray*}
The proof is now complete. \qed

\vspace{.1in}

\section{Martin Boundary}

Superharmonic and harmonic functions with respect to $X^D$ 
have been studied in the context of general
theory of Markov processes and their potential theory
(see, for instance, Kunita-Watanabe \cite{wa:kun}). 
From the
general theory, we know that
positive harmonic functions with respect to $X^D$ 
admit Martin representations. However, 
no particular attention was 
paid to the special case of
harmonic functions with respect to the killed stable process.
For instance,
the relationship between the Martin boundary of $X^D$
and the Euclidean boundary $\partial D$ of $D$ has not been studied.
 
In this section we assume that $D$ is a bounded Lipschitz domain.
In the first part of this section we 
are going to show that the Martin
boundary of $X^D$ and the Euclidean boundary $\partial D$ coincide. 
Our proof of the identification
between the Martin boundary and the Euclidean boundary is 
similar to the argument of Bass--Burdzy \cite{bu:bas}
in the Brownian motion case.

Fix $x_0\in D$ and set
\[
M_D(x, y)=\frac{G_D(x, y)}{G_D(x_0, y)},  \qquad   x, y\in D.
\]
The Martin boundary is the set $\partial _MD=D^*\setminus D$,
where $D^*$ is the smallest compact set for which $M_D(x, y)$ is
continuous in $y$ in the extended sense. 

\begin{lemma}\label{2A}
Suppose that $D$ is a bounded Lipschitz domain.
Then for any $\gamma>0$, 
\[
\lim_{x\to\partial D}G_D(x, y)=0
\]
uniformly on $D_{\gamma}=\{y\in D: \delta(y, \partial D)\ge\gamma\}$.
\end{lemma}

\pf
Suppose that 
${\cal C} =\{x=(x_1, \dots, x_n)\in \RR^n: 0<x_n,
(x^2_1+\cdots+x^2_{n-1})^
{1/2}<\eta x_n\}$ is a cone with vertex at the origin $O$. 
For any $r>0$, set
\[
{\cal C}_r={\cal C}\cap B(O, r)
\]
and
\[
T_{{\cal C}_r}=\inf\{t>0, X_t\in\overline{{\cal C}_r}\}.
\]
One can easily show (similar to the proof of Proposition 1.19 of
\cite{kl:chu}) that for any $t>0$, the function
\[
x\mapsto P_x(t<T_{{\cal C}_r})
\]
is upper semi--continuous in $\RR^n$. Thus for any $s>0$,
we have
\[
\limsup_{x\to O} P_x(T_{{\cal C}_r}>s)\le P^O(T_{{\cal C}_r}>s)=0
\]
since $P^O(T_{{\cal C}_r}=0)=1$. Now use the fact that $D$
satisfies the uniform exterior cone condition we can easily
see that
\[
\lim_{x\to z\in \partial D}P_x(\tau_D>s)=0
\]
uniformly in $z\in\partial D$, i.\/e., for any $\epsilon>0$
there exists $\delta'>0$
such that 
\begin{equation}\label{eqn:21}
P_x(\tau_D>s)<\epsilon, \qquad \mbox{ if } \delta(x, \partial D)<\delta'.
\end{equation}

We know that
\[
G_D(x, y)=G(x, y)-E_x[G(X_{\tau_D}, y)]
\]
and that $G(x, y)$ is bounded $D^c\times D_{\gamma}$. 
Now use the fact (\ref{eqn:21}) and argue along the line
of the proof of Theorem 1.23 of \cite{kl:chu} we
easily arrive at our conclusion.
\qed

\vspace{.1in}

Take a positive number $\epsilon<\delta(x_0, \partial D)/4$.

\begin{lemma}\label{2a}
Suppose $x\in D$ with $|x-x_0|>4\epsilon$. There exists a 
constant $c_1=c_1(\epsilon, D, x, x_0)$ such that 
\[
M_D(x, y)\le c_1 \quad \mbox{ for } y\in
D\setminus \left( \overline{B(x_0, \epsilon)}\cup\overline{B(x, \epsilon)}
 \right).
\]
\end{lemma}

\pf
Pick $y_0\in \partial B(x_0, 2\epsilon)$. By the explicit
formula for the Green function of balls (see \cite{g1:blu} for
instance) we know that
\[
G_D(x_0, y_0)\ge G_{B(x_0, 3\epsilon)}(x_0, y_0)\ge \delta(\epsilon)>0.
\]

On the other hand, we know that
$G_D(x, y_0)\le 2^{-\alpha}\pi^{-n/2}\Gamma(
(n-\alpha)/2)\Gamma(\alpha/2)^{-1}|x-y_0|^{\alpha-n}$.
 Therefore
$G_D(x, y_0)$ is bounded above by a constant depending on $\epsilon$
in $x\in D$ with $|x-x_0|>4 \epsilon $.
Thus $M_D(x, y_0)$ is bounded above in 
 $x\in D$ with $|x-x_0|>4 \epsilon $. But from the boundary Harnack
principle (see \cite{kr:bog}) we get that $M_D(x, y)$ is comparable
to $M_D(x, y_0)$ for all points $y$ in $D\setminus(B(x, \epsilon)
\cup B(x_0, \epsilon))$. The lemma follows.
\qed

\vspace{.1in}

\begin{lemma}\label{2b}
Let $x$, $x_0$, $\epsilon$ be as above. Then $M_D(x, y)$ is
a H\"older continuous function of $y$ in
$D\setminus \left(\overline{B(x_0, \epsilon)}\cup\overline{B(x, \epsilon)}
\right)$ with 
 H\"{o}lder exponent and coefficient
 depending only on $x$, $x_0$, $\epsilon$ and $D$.
\end{lemma}

\pf
For any set $A$, we define that
\[
\mbox{Osc}_A\,f=\sup_{y\in A} f(y)-\inf_{y\in A}f(y).
\]
Let $f(y)=M_D(x, y)$. Let $y_0\in D_{\epsilon}=
D\setminus \left(\overline{B(x_0, \epsilon)}\cup\overline{B(x, \epsilon)}
\right)$.
Since by Lemma \ref{2a} $f$ is bounded by $c_1$ on $D_{\epsilon}$,
 $\mbox{Osc}_{D_{\epsilon}}\, f \, \le c_1$. So
it suffices to show that there exists $\rho=\rho(\epsilon, D, x, x_0)<1$
such that
\begin{equation}
\mbox{Osc}_{D\cap B(y_0, r)}\,f\le \rho \, \mbox{Osc}_
{D\cap B(y_0, 2r)}\, f
\quad \mbox{ for }r<\epsilon/4.
\end{equation}

Suppose $r<\epsilon/4$, and let $g$ be the ratio of any two positive harmonic
functions on $D_{\epsilon/4}$ vanishing continuously on
$D^c$. By considering $ag+b$ for suitable $a$ and $b$, we may assume 
\[
\sup_{D\cap B(y_0, 2r)}g=1,\,\,\,\, \inf_{D\cap B(y_0, 2r)}g=0.
\]

If $\sup_{D\cap B(y_0, r)}g\le 1/2$, then since $\inf_{D\cap 
B(y_0, 2r)}g\ge 0$, we have
\[
\mbox{Osc}_{D\cap B(y_0, r)}\,g\le\frac12
=\frac12 \mbox{Osc}_{D\cap B(y_0, 2r)}\,g .
\]
If $\sup_{D\cap B(y_0, r)}g\ge 1/2$, there exists a point $y_1$ in 
$D\cap B(y_0, r)$ with $g(y_1)\ge 1/2$. But then by the boundary Harnack
principle with $V=\{x: \delta(x, D)<\epsilon\}\setminus(\overline{B}
(x_0, r)\cup\overline{B}(x, r))$ and $K=\overline{D}\setminus(\overline{B}
(x_0, r)\cup\overline{B}(x, r))$
, 
there exists a constant $c_2=c_2(\epsilon, D, x, x_0)\in (0, 1)$
such that
\[
\inf_{D\cap B(y_0, r)}g\ge c_2g(y_1).
\]
Since $\sup_{D\cap B(y_0, r)}g\le 1$, in this case we have
\[
\mbox{Osc}_{D\cap B(y_0, r)}\,g\le 1-\frac{c_2}2 =
\left( 1-\frac{c_2}2 \right) \mbox{Osc}_{D\cap B(y_0, 2r)}\, g .
\]
So we have (3.1) with
$\rho=\max \left\{\frac12, 1-\frac{c_2}2 \right\} $.
Therefore $M_D(x, y)$ is
a (globlly) H\"older continuous in $y \in
D\setminus \left(\overline{B(x_0, \epsilon)}\cup\overline{B(x, \epsilon)}
\right)$.
\qed

\vspace{.1in}

A direct consequence of Lemma \ref{2b} is that 
$M_D(x, y)=G_D(x, y)/G_D(x_0, y)$
converges when  $y\to z\in \partial D$. 
Let the limit be denoted as $M_D(x, z)$. This implies that
the Martin boundary of $D$
can  identified with a subset of $\partial D$.

It is also well known that for a bounded Lipschitz domain $D$
with Lipschitz characteristic constants
$(r_0, \, A_0)$,
there exists $\kappa=\kappa(A_0)\in (0, 1)$ such that for every
$\epsilon \in(0, r_0)$ and $z\in\partial D$, there is a point
$A_\epsilon (z)\in D\cap B(z, r)$ such that $B(A_\epsilon (z), \kappa r)
\subset D\cap B(z, r)$. It is not difficult to show the following
(cf. Lemma 6 of Bogdan \cite{bog2}).

\begin{lemma}\label{2B}
For any $z\in \partial D$, $M_D(\cdot, z)$ is harmonic
with respect to $X^D$.
\end{lemma}

\pf Clearly  any fixed $x\in D$ and and $r<\delta(x, \partial D)$, 
\[ 
M_D(x, y)= E_x \left[ M_D(X_{\tau_{B(x, r)}}, y); \, \tau_{B(x, r)}
< \tau_D  \right]  \quad \hbox{for } y\in D\setminus \overline{B(x, r)}.
\]
In particular, 
\begin{equation}
M_D(x, A_{\epsilon}(z)) =
E_x \left[ M_D(X_{\tau_{B(x, r)}},  A_{\epsilon}(z)); \, \tau_{B(x, r)}
< \tau_D  \right]\label{BB} 
\end{equation}
for any $0<\epsilon <\min\{r, \, r_0 \}$.
By Fatou's lemma, 
$$
M_D(x, z)\ge E_x \left[ M_D(X_{\tau_{B(x, r)}}, z); \, \tau_{B(x, r)}
< \tau_D  \right]. 
$$
Therefore $M_D(X_{\tau_{B(x, r)}}, z)$ is $P_x$-integrable. 
Put
\[
\epsilon_0=\min\left\{\frac{\delta(x_0, \partial D)}4, \frac{r_0}2,
\frac{r}4\right\},
\]
Then by Lemma 13 of \cite{kr:bog} we get that there exists
$C_1=C_1(D)>0$ such that for any $y\in D\cap
B(z, \epsilon_0)$ and $\epsilon\in(0, \epsilon_0)$,
\[
M_D(w, A_{\epsilon}(z))\le C_1 M_D(w, y) \quad \hbox{for } 
w\in D\setminus B(z, 2\epsilon).
\]
Letting $y\to z$ we get that for any $\epsilon\in(0, \epsilon_0)$,
\begin{equation}
M_D(w, A_{\epsilon}(z))\le C_1 M_D(w, z) \quad \hbox{for } 
w\in D\setminus B(z, 2\epsilon).\label{1}
\end{equation}

For $w\in D\cap B(z, 2\epsilon)$, $|w-x|>3r/2$ and thus by the
explicit formula for $K_{B(x, r)}$ we know that there is a constant
$C_2=C_2(r)>0$ such that
$K_{B(x, r)}(x, w)\le C_2$ for $w\in D\cap B(z, 2\epsilon)$. 
Hence for any $w\in D\cap B(z, 2\epsilon)$,
\begin{eqnarray*}
&&E_x \left[ M_D(X_{B(x, r)}, \,  A_{\epsilon}(z)); \, 
X_{B(x, r)} \in D\cap B(z, 2\epsilon) \right] \\
&\le&\frac{C_2}{G_D(x_0, A_{\epsilon}(z))}\int_{D\cap B(z, 2\epsilon)}
G(w, A_{\epsilon}(z))dw\\
&\le&\frac{C_3}{G_D(x_0, A_{\epsilon}(z))}\int_{D\cap B(z, 2\epsilon)}
|w-A_{\epsilon}(z)|^{\alpha-n}dw\\
&\le&\frac{C_4}{G_D(x_0, A_{\epsilon}(z))}\epsilon^{\alpha}.
\end{eqnarray*}
From Lemma 5 of \cite{kr:bog} we know that there
exists a constant $C_5=C_5(D, x_0)>0$ and positive
number $\gamma=\gamma(D)<\alpha$ such that
$
G_D(x_0, A_{\epsilon}(z))\ge C_5 \epsilon^{\gamma} $.
Therefore
\begin{equation}
E_x \left[ M_D(X_{B(x, r)}, \,  A_{\epsilon}(z)); \, 
X_{B(x, r)} \in D\cap B(z, 2\epsilon) \right]
\le C_4 C_5^{-1}\epsilon^{\alpha-\gamma}.\label{2}
\end{equation}
Now combine (\ref{1}) and (\ref{2}) we see that the family of functions
$\{M_D(X_{B(x, r)}), \, A_{\epsilon}(z)
): \, 0<\epsilon < \epsilon_0 \}$ is $P_x$-uniformly integrable. Letting $\epsilon\to 0$ in 
(\ref{BB}) yields
\[
M_D(x, z)= E_x \left[ M_D(X_{\tau_{B(x, r)}}, z); \, \tau_{B(x, r)}
< \tau_D  \right] \quad \hbox{ for any } r < \delta ( x , \partial D). 
\]
Thus $M_D(\cdot, z)$
is harmonic with respect to $X^D$.
\qed

\vspace{.1in}

The next result tells that
each Euclidean boundary point corresponds to a different non-negative
harmonic function. Hence the Martin boundary can not be identified
with a proper subset of the Euclidean boundary.

\begin{lemma}\label{2c}
If $M_D(\cdot, z_1) 
\equiv M_D(\cdot, z_2)$ for $z_1$, $z_2\in \partial D$,
then $z_1=z_2$.
\end{lemma}

\pf
Let $\epsilon>0$ be such that
\[
\epsilon<\min\left\{ \frac{\delta(x_0, \, \partial D)}4, \, 
\frac{r_0}2\right\},
\]
where $r_0$ comes
from the Lipschitz characteristic constants $(r_0, A_0)$
of $D$.
First we are going to
show that $M_D(x, w)\to 0$
uniformly in $w\in\partial D$ 
as $\delta(x, \partial D\setminus B(w, 3\epsilon))\to 0$. 
In fact,  for any given  $\eta>0$,
by Lemma \ref{2A} there is a $\beta=\beta (\eta, A_0, \epsilon)>0$ such that
$G_D(x, A_{\epsilon}(w))<\eta$ for $ w\in\partial D$ and 
$x\in D$ with  $\delta(x, \partial D)< \beta$. 
Let $D_0$ be a  smooth domain such that
\[
\left\{x\in D; \delta(x, \partial D)\ge\frac{\kappa \epsilon}2 \right\}\subset
D_0\subset\overline{D_0}\subset D.
\]
Then for all $w\in \partial D$,
$\delta(A_{\epsilon}(w), \partial D_0)>\kappa\epsilon/2$ and so
by Theorem 1.1 of \cite{s1:che}  
\begin{eqnarray*}
G_D(x_0, A_{\epsilon}(w))&\ge& G_{D_0}(x_0, A_{\epsilon}(w))\\
&\ge&C(D_0)\min\left\{\frac1{|x_0-A_{\epsilon}(w)|^{n-\alpha}}, \, \frac{
\delta(x_0, \partial D_0)^{\alpha/2}\delta(A_{\epsilon}(w), 
\partial D_0)^{\alpha/2}}{|x_0-A_{\epsilon}(w)|^{n}}\right\}\\
&\ge&C(D_0)\min\left\{\frac1{d_D^{n-\alpha}}, \, \frac{(3\epsilon)^{\alpha/2}
(\kappa\epsilon/2)^{\alpha/2}}{d_D^n}\right\}\\
& := &C_1>0,
\end{eqnarray*}
where $d_D$ is the diameter of $D$. Therefore 
\[
M_D(x, A_{\epsilon}(w))<\eta/C_1, \qquad \forall w\in\partial D,
\]
whenever $\delta(x, \partial D)< \beta$. 
Fix $w\in\partial D$. Clearly $|x-w|>2\epsilon$
for any $x\in D$ with $\delta(x, \partial D\setminus B(w, 3\epsilon))
<\epsilon$. Now by Lemma
13 of \cite{kr:bog} with $r=\epsilon$, we get that
there is a $C_2=C_2(A_0)>0$ such that
\[
M_D(x, y)\le C_2M_D(x, A_{\epsilon}(w)) \quad \hbox{for } y\in D\cap 
B(w, \epsilon)
\]
whenever $x\in D$ satisfies $\delta(x, \partial D\setminus B(w, 3\epsilon))
<\epsilon$.
Therefore
\[
M_D(x, y)\le \frac{C_2}{C_1}\eta \quad \hbox{for } y\in D\cap 
B(w, \epsilon)
\]
whenever $x\in D$ satisfies $\delta(x, \partial D\setminus B(w, 3\epsilon))
<\epsilon\wedge\beta$.
Letting $y\to w$ we get that
\[
M_D(x, w)\le \frac{C_2}{C_1}\eta
\]
for $x\in D$ with $\delta(x, \partial D\setminus B(w, 3\epsilon))
<\epsilon\wedge\beta$.
Therefore
\begin{equation}\label{eqn:unif}
M_D(x, w)\to 0 \hbox{ 
uniformly in } w\in\partial D \hbox{ as }
\delta(x, \partial D\setminus B(w, 3\epsilon))\to 0.
\end{equation}

Suppose that $M_D(\cdot, w)=M_D(\cdot, z)$ for some $w$, $z\in \partial D$, 
$w\neq z$, and let $\epsilon<|w-z|/8$. By the above argument, 
$M_D(x, w)\to 0$ uniformly when 
$\delta(x, \partial D\setminus B(w, 2\epsilon))\to 0$
or when $\delta(x, \partial D\setminus B(z, 2\epsilon))\to 0$.
Therefore $M_D(x, w)\to 0$ uniformly as $\delta(x, \partial D)\to 0$.
Since $M_D(\cdot, w)$ is a non-negative harmonic function with
respect to $X^D$ which continuously vanishes on $\partial D$,
it must be identically zero by Theorem \ref{1e}. This contradicts the
fact that $M_D(x_0, w)=1$. The proof is now complete.
\qed

\vspace{.1in}

Combining the lemmas above we get the following result.

\begin{thm}\label{2d}
The Martin boundary of $D$ can be identified with its Euclidean boundary
$\partial D$.
\end{thm}

\begin{thm}\label{2e}
For each $z\in \partial D$, $M_D(x, z)$ is minimal harmonic with respect
to $X^D$.
\end{thm}

\pf
Fix $z\in\partial D$ and suppose $h\le M_D(\cdot, z)$, where
$h$ is a positive harmonic function with respect to $X^D$. 
By Theorem \ref{2d} we know that there is a measure $\mu$ on 
$\partial D$ such that
\[
h(\cdot)=\int_{\partial D}M_D(\cdot, w)\mu(dw).
\]
If $\mu$ is not a multiple of the point mass at $z$, then there
is a finite measure $\nu\le \mu$ such that $\delta(z, \mbox{supp}(\nu))
>0$. Let
\[
u(\cdot)=\int_{\partial D}M_D(\cdot, w)\nu(dw).
\]
Then $u$ is a positive harmonic function with respect to $X^D$ bounded
by $M_D(\cdot, z)$.

Recall from (\ref{eqn:unif} in 
the proof of Lemma \ref{2c} that $M_D(x, z) \to 0$ uniformly as
$\delta(x, \partial D\setminus B(z, \epsilon) \to 0$. So the same
is true of $u$. But for each $w\in\mbox{supp}(\nu)$, we can see that
$M_D(x, w)\to 0$ uniformly as $\delta(x, \partial D\cap B(z, 2\epsilon)
\to 0$ provided $2\epsilon<\delta(z, \mbox{supp}(\nu))$. So it follows by
the dominated convergence theorem that $u(x)\to 0$ as 
$\delta(x,\partial D\cap B(z, 2\epsilon)
\to 0$. But then $u$ is a positive harmonic function of $X^D$
which continuously vanishes on $\partial D$. This implies that
$\nu$ is 0, or that $\mu=c\delta_z$ for some $c$. 
\qed

\vspace{.1in}

From Theorem \ref{2d} and the general theory of Martin representation
(cf. \cite{wa:kun}), we have

\begin{thm}\label{2f}
If $D$ is a bounded 
Lipschitz domain in $\RR^n$, then the restriction to $D$ of 
any positive superharmonic 
function $f$ with respect to $X^D$ can be written uniquely as
\begin{equation}\label{eqn:martin}
f(x)=\int_DG_D(x, y)\nu(dy)+\int_{\partial D}M_D(x, z)\mu(dz), 
\end{equation}
where $\nu$ and $\mu$ are finite measures on $D$ and $\partial D$
respectively.
\end{thm}

When $D$ is a bounded $C^{1, 1}$ domain, we can say more
about the Martin kernel of $D$.

\begin{thm}\label{2g}
Suppose that $D$ is a bounded $C^{1, 1}$ domain, then $M_D(\cdot, \cdot)$
is continuous function on $D\times \partial D$. Furthermore, there
exists a constant $C=C(D, \alpha)>0$ such that for any $x\in D$ and $z\in
\partial D$,
\[
\frac{\delta(x, \partial D)^{\alpha/2}}{C|x-z|^n}\le M_D(x, z)\le
\frac{C\delta(x, \partial D)^{\alpha/2}}{|x-z|^n}.
\]
\end{thm}

\pf
The joint continuity of $M_D$ follows from the definition
of Martin kernel and Lemma
\ref{2b}. The estimates on $M_D$ follows easily from
Theorems 1.1 and 1.2 of Chen and Song \cite{s1:che}.
\qed

\vspace{.1in}

From Theorems 1.1--1.2 of Chen and Song 
\cite{s1:che} and Theorem \ref{2g} above
we have

\begin{thm}\label{2n} (3G Theorem).
Suppose that $D$ is a bounded $C^{1, 1}$ domain, then there
exists $C=C(D, \alpha)>0$ such that
\[
\frac{G_D(x, y)M_D(y, z)}{M_D(x, z)}\le C\frac{|x-z|^{n-\alpha}}
{|x-y|^{n-\alpha}|y-z|^{n-\alpha}}, \qquad x, y\in D, z\in\partial D.
\]
\end{thm}

\pf The proof is the same as the proof of Theorem 1.6 of 
Chen and Song \cite{s1:che}
and we omit it here.
\qed

\vspace{.1in}

Using Theorems \ref{2g} and \ref{2n}, we can prove a conditional gauge 
theorem, which complements the two conditional gauge theorems established in 
Chen and Song 
\cite{s2:che}. Before we state and prove the conditional gauge theorem,
we need to do some preparations first.

\begin{defn}\label{def:kato}
A Borel measurable function $q$ on $\RR^n$ is said to be in the Kato class 
$\KK_{n, \alpha}$ if
\begin{equation}\label{eqn:kato}
\lim_{r\downarrow 0}\sup_{x\in \RR^n}\int_{|x-y|\le r}
\frac {|q(y)|}{|x-y|^{n-\alpha}}dy=0.
\end{equation}
\end{defn}

\smallskip

For $q\in\KK_{n, \alpha}$, set
\[
e_q(t)=\exp \left(\int^t_0q(X_s)ds \right).
\] 
From Chen and Song \cite{s2:che}, we know that the following semigroup
\[
T_tf(x)=E_x[e_q(t)f(X_t); \, t<\tau_D], \qquad x\in D,
\]
admits an integral kernel $k_q(t, x, y)$. The function
\[
g(x) := E_x[e_q(\tau_D)]
\]
is called the gauge function of $(D, q)$. 
It is shown in \cite{r2:chu} that  either
$g$ is  identically identically
infinite or $g$ is bounded on $D$. In the latter case, $(D, q)$
is said to be {\it gaugeable}.
When $(D, q)$ is gaugeable,
it can be showb (see \cite{s2:che}) that
\[
V_q(x, y)=\int^{\infty}_0k_q(t, x, y)dt, \qquad x, y\in D,
\]
is well defined and is continuous off the diagonal.

Suppose that $h>0$ is a positive superharmonic function with respect to
$X^D$. Note that by Theorem \ref{1c} above, we have 
(see,  e.g., page 11 of Dynkin \cite{eb:dyn}) that
$$ h(x) \geq E^x [h(X^D_t)].
$$
We define
\[
p_D^h(t, x, y)=h(x)^{-1}p_D(t, x, y)h(y), \qquad t>0, x, y\in D,
\]
where $p_D$ is the transition density function of killed symmetric stable
process $X^D$ in $D$. It is easy to check that $p_D^h$ is a transition
density and it determines a Markov process on the state space
$D_{\partial}=D\cup \{\partial\}$. This process is called the 
$h$--conditioned symmetric stable process. Similar to Propositions
5.2--5.4 of Chung and Zhao \cite{zh:chu}, we have the following

\begin{lemma}\label{2i}
For $t>0$, if $\Phi\ge 0$ is an ${\cal F}_t$--measurable function, then
\[
E_x^h[\Phi; t<\tau_D]=h(x)^{-1}E_x[\Phi\cdot h(X_t); t<\tau_D], 
\qquad x\in D.
\]
\end{lemma}

Recall that  $\{ {\cal F}_t , \, t \geq 0\}$ be the minimal admissible
$\sigma$-fields generated by $X$.
For any stopping time $T$ of  $\{ {\cal F}_t , \, t \geq 0\}$, 
${\cal F}_{T+}$ is the class
of subsets $\Lambda$ of ${\cal F}$ such that
\[
\Lambda\cap\{T\le t\}\in{\cal F}_{t+}, \qquad t\ge0.
\]
${\cal F}_{T-}$ is the $\sigma$--field generated by ${\cal F}_{0+}$ and
the class of sets
\[
\{t<T\}\cap \Lambda, \qquad t\ge0, \lambda\in{\cal F}_t.
\]

\begin{lemma}\label{2j}
For any stopping time $T$ and any ${\cal F}_{T+}$--measurable
function $\Phi\ge0$, 
\[
E_x^h[\Phi; T<\tau_D]=h(x)^{-1}E_x[\Phi\cdot h(X_T); T<\tau_D].
\]
\end{lemma}

\begin{lemma}\label{2k}
For any stopping time $T$, $A\in {\cal F}_{T+}$ and any ${\cal F}_{
\tau_{D-}}$--measurable variable $\Phi\ge 0$,
\[
E_x^h[A\cap(T<\tau_D); \Phi \circ \theta_T]=E_x^h[A\cap(T<\tau_D);
E^{X_T}_h(\Phi)],
\]
where $\theta_t $ is the shift operator for process $X$.
\end{lemma}

Now let $D$ be a bounded Lipschitz domain. For each $z\in\partial D$,
the $M_D(\cdot, z)$--conditioned symmetric stable process will be
called the $z$-symmetric stable process, and the associated probability
and expectation  will be denoted by $P_x^z$ and $E_x^z$, respectively. 

For any $y\in D$, $G_D(\cdot, y)$ is harmonic in $D\setminus\{y\}$
with respect to $X^{D\setminus\{y\}}$. Hence we can define the
$G_D(\cdot, y)$--conditioned symmetric stable process on the state
space $(D\setminus\{y\})\cup \{\partial\}$, with lifetime
$\tau_{D\setminus\{y\}}$. It will be referred to as the $y$--conditioned
symmetric stable process, and the associated probability and expectation
will be denoted by $P_x^y$ and $E_x^y$ respectively.
The following result immediately follows from Theorem \ref{2n}.

\begin{corollary}\label{cond} (Conditional Lifetime)
Suppose that $D$ is a bounded $C^{1, 1}$ domain. Then
$$ \sup_{x\in D, \, z\in \partial D} E_x^z [ \tau_D ] < \infty .
$$
\end{corollary}

\begin{thm}\label{2o} (Conditional Gauge Theorem).
Suppose that $D$ is a bounded $C^{1, 1}$ domain and $q\in\KK_{n, \alpha}$.
If $(D, q)$ is gaugeable,
then there exists $c>1$ such that
\[
c^{-1}\le\inf_{x\in D, z\in\partial D}E_x^z[e_q(\tau_D)]
\le \sup_{x\in D, z\in\partial D}E_x^z[e_q(\tau_D)]\le c.
\]
\end{thm}

\pf Suppose $x, y\in D$ and $z\in\partial D$. For any $w\in D$, by
it follows from Lemma 6.5 of Chen and Song \cite{s2:che} that
\begin{eqnarray*}
\lim_{y\to z}\frac1{G_D(x, y)}V_q(x, w)q(w)G_D(w, y)
&=&\lim_{y\to z}V_q(x, w)q(w) \frac{G_D(w, y)/
G_D(x_0, y)}
{G_D(x, y)/G_D(x_0, y)}\\
&=&V_q(x, w)q(w) \frac{M_D(w, z)}{M_D(x, z)}.
\end{eqnarray*}
Now from Theorem 1.6 (3G Theorem) of Chen and Song 
\cite{s1:che} and Theorem 5.2
of Chen and Song \cite{s2:che} we have that
\[
\{V_q(x, \cdot)G_D(\cdot, y)|q(\cdot)|/G_D(x, y): x, y\in D\}
\]
is uniformly integrable. Hence it follows from Theorem 5.4
of Chen and Song \cite{s2:che} that
\[
\lim_{y\to z}E_x^y[e_q(\tau_{D\setminus\{y\}})]=
1+\frac1{M_D(x, z)}\int_DV_q(x, u)q(u)M_D(u, z)du.
\]
However, one can show, by using
an argument similar to the proof of  
Theorem 5.4 of Chen and Song \cite{s2:che}, that
\[
E_x^z[e_q(\tau_D)]=1+\frac1{M_D(x, z)}\int_DV_q(x, u)q(u)M_D(u, z)du.
\]
Therefore
\begin{equation}
\lim_{y\to z}E_x^y[e_q(\tau_{D\setminus\{y\}})]=
E_x^z[e_q(\tau_D)]
\end{equation}
 and the theorem now follows from Theorem 5.6 of Chen and Song 
 \cite{s2:che}.
\qed.

\medskip

\begin{thm}\label{2p}
Suppose that $D$ is a bounded $C^{1, 1}$ domain and $q\in\KK_{n, \alpha}$.
If $(D, q)$ is gaugeable,
then for any fixed point $x_0 \in D$ and $z\in \partial D$, 
\begin{equation}\label{eqn:3.3}
\lim_{D\ni y\to z}\frac{V_q(x, y)}{V_q(x_0, y)}=
\frac{E_x^z[e_q(\tau_D)]}{E^{x_0}_z[e_q(\tau_D)]} \, M_D(x, z).
\end{equation}
Furthermore, 
\begin{equation}\label{eqn:3.4}
\lim_{D\ni y\to z}\frac{V_q(x, y)}{\delta(y, \partial D)
^{\alpha/2}}=E_x^z[e_q(\tau_D)] \, 
\lim_{D\ni y\to z\in\partial D}\frac{G_D(x, y)}{\delta(y, \partial D)
^{\alpha/2}}.
\end{equation}
\end{thm}

\pf Since  
\begin{equation}\label{eqn:3.5}
E_x^y[e_q(\tau_{D\setminus\{y\}})]=V_q(x, y)/
G_D(x, y) 
\end{equation}
for $x, y \in D$ by Theorem 5.5 of Chen and Song \cite{s2:che}
and so (\ref{eqn:3.3}) follows immediately from it.
Identity (\ref{eqn:3.4}) follows from 
(\ref{eqn:3.5}) and from Lemma 6.5 of Chen and Song 
\cite{s1:che} which asserts that the limit
$$
\lim_{D\ni y\to z\in\partial D}\frac{G_D(x, y)}{\delta(y, \partial D)
^{\alpha/2}}
$$
exists and forms a positive and continuous function in
$(x, z)\in D\times \partial D$.
 \qed

\vspace{.1in}

As a consequence of Theorems \ref{2o} and \ref{2p} we get that,
for a bounded $C^{1, 1}$ domain $D$ and a $q\in\KK_{n, \alpha}$, 
if $(D, q)$ is gaugeable, then the
Martin kernel of the generalized Schr\"odinger
operator $-(-\Delta)^{\alpha/2}+q$ with zero exterior
condition on $D^c$ is comparable to
Martin kernel $M_D$.

\medskip

\noindent{\bf Remark 3.1} 
Recently in \cite{s3:che}
we were able to extend the 3G theorem and 
conditional gauge theorem established for bounded $C^{1,1}$ domains
in Chen and Song \cite{s1:che}, \cite{s2:che} to bounded Lipschitz domains.
Thus Theorems \ref{2n} and \ref{2o} and (\ref{eqn:3.3})
in Theorem \ref{2p} in fact hold on bounded Lipschitz
domains as well. Thus under the condition that $D$ is a bounded 
Lipschitz domain and 
$(D, q)$ is gaugeable, 
the Martin kernel of the generalized Schr\"odinger
operator $L= -(-\Delta)^{\alpha/2}+q$ with zero exterior
condition on $D^c$ is comparable to
Martin kernel $M_D$ and the Martin boundary for $L$
coincides with the Euclidean boundary $\partial D$ of $D$.

\medskip

When $D$ is a ball, we can actually get an explicit formula for the
Martin kernel of $D$. This follows easily form the definition
of the Martin kernel and the explicit formula for the Green functions
of balls (see Corollary 4 of Blumenthal,  Getoor and  Ray
\cite{g1:blu}). We record this fact as follows.

\vspace{.1in}

\noindent {\bf Example}. 
If $B=B(0, r)$, then
\[
M_B(x, w)=\frac{(r^2-|x|^2)^{\alpha/2}}{|x-w|^n}, \qquad x\in D, w\in
\partial D.
\]

\vspace{.1in}

From the formula above, we know that
\[
h(x)=\int_{\partial B}\frac{(r^2-|x|^2)^{\alpha/2}}{|x-w|^n}dw, \qquad x\in D
\]
is a positive harmonic function with respect to $X^B$. From Theorem \ref{1e}
we know that $h$ can not be a bounded function on $B$.
In fact one can check directly in this case that for each $z\in\partial B$ 
\[
\lim_{B\ni x\to z }h(x)=\infty.
\]

In the Brownian motion case, the Martin boundary can be approached along
Brownian paths. While for a symmetric stable process, we 
know from Lemma 6 of Bogdan \cite{kr:bog}
that, with probability 1, it will not hit $\partial D$
upon first exiting from a bounded domain $D$ satisfying
the uniform exterior cone condition. Our next 
theorem gives the relationship between the Martin boundary
and the (conditioned) stable paths.

\begin{thm}\label{2l}
Suppose $D$ is a bounded Lipschitz domain. Then for every $x\in D$
and $z\in\partial D$,
\begin{eqnarray*}
P_x^z\{\tau_D<\infty\}&=&1;\\
P_x^z\{\lim_{t\uparrow\tau_D}X(t)=z\}&=&1.
\end{eqnarray*}
\end{thm}

\pf
Let $z\in\partial D$, $r_m\downarrow 0$, $B_m=(z, r_m), D_m=D\setminus
\overline{B_m}$ and set
\[
T_m=\inf\{t>0: X_t\in B_m\}, \qquad R_m=\tau_{B_m\cap D}.
\]

We may assume that $x\in D_m$ for $n\ge1$. By Lemma \ref{2b}, $M_D(\cdot, z)$
can be continuously extended onto $\overline{D_m}$ by setting $M_D(w, z)=0$
for $w\in \partial D\setminus B_m$. Since $M_D(\cdot, z)$ is harmonic in $D_m$
with respect to $X^D$, we have by Theorem \ref{1d} and Lemma \ref{2j}
\begin{eqnarray*}
M_D(x, z)&=&E_x\left[ M_D(X_{\tau_{D_m}}, z)\right] \\
&=&
E_x\left[ M_D(X_{T_m}, z); T_m<\tau_D\right] \\
&=& M_D(x, z)P_x^z\left( T_m<\tau_D\right) .
\end{eqnarray*}
It follows that for all $n\ge1$ we have
\begin{equation}\label{eqn:3.7}
P_x^z\{T_m<\tau_D\}=1.
\end{equation}

Note that for each fixed $z\in \partial D$,
$M_D(x, z)$ is bounded in $x\in B^c_k\cap D$ by continuity.
Let $C_k$ denote its bound. Applying Lemmas \ref{2j} and \ref{2k} 
twice, we have  for all $k<m$:
\begin{eqnarray}
P_x^z\{T_m<\tau_D,\,  R_k\circ\theta_{T_m}<\tau_D\}
&=&E_x^z\left[ P^{X_{T_m} }_z[R_k<\tau_D]; \, T_m<\tau_D \right]\nonumber\\
&=&\frac1{M_D(x, z)}E_x\left[ M(X_{T_m}, z)
P^{X_{T_m}}_z[R_k<\tau_D]; \, T_m<\tau_D \right] \nonumber\\
&=&\frac1{M_D(x, z)}E_x\left[ E^{X_{T_m} }[
M_D(X_{R_k}, z); \, R_k<\tau_D]; \, T_m<\tau_D\right]\nonumber\\
&\le& \frac{C_k}{M_D(x, z)}P_x\left( T_m<\tau_D\right)\label{eqn:3.8} .
\end{eqnarray}
By the definition of $T_m$ and the
quasi left continuity of the unconditioned process $X$, we have
\begin{eqnarray*}
\lim_{m\to\infty}P_x\{T_m<\tau_D\} &\le& P_x\{
\lim_{m\to\infty}T_m\le \tau_D\}\\
&\le& P_x\{T_{\{z\}}\le \tau_D\} =0
\end{eqnarray*}
because  $z\in \partial D$ and by (\ref{eqn:2.1}) with $D$ in place of $U$
that $P_x( X_{\tau_D} \in \partial D) =0$.
It follows from (\ref{eqn:3.8}) that the
left hand side there converges to zero as $m\to\infty$ for each $k$.
Therefore there exists a subsequence $\{m_j\}$ such that
\[
\sum_{j=1}^{\infty}P_x^z\{T_{m_j}<\tau_D; R_k\circ\theta_{T_{m_j}}
<\tau_D\}<\infty,
\]
and consequently by Borel--Cantelli lemma we have
\[
P_x^z\{[ T_{m_j}<\tau_D; R_k\circ\theta_{T_{m_j}}
<\tau_D] \ \mbox{ infinitely often}\, \}=0.
\]
Together with (\ref{eqn:3.7}) this implies that for $k\geq 1$ and $P_x^z$--a.\/e. $\omega$, 
there exists an integer $N(\omega)<\infty$ such that
\[
X_t(\omega)\in B(z, r_k) \mbox{ for all } t\in[T_{N(\omega)}(\omega),
\tau_D(\omega)).
\]
For each $k$ let $N(k)$ be the smallest $N$ for which the above
is true. Then $T_{N(k)}\uparrow\tau_D$; otherwise, we would have
$X_t=z$ for all $t\in[\lim_{k\to\infty}T_{N(k)}, \tau_D)$,
which is impossible since $z\notin D_{\partial}$. This proves that
$X_t\to z$ as $t \uparrow \tau_D$.
\qed

\vspace{.1in}

Functions harmonic in $D$ with respect to $X^D$ do not come
from solving Dirichlet exterior problems. Therefore the usual probabilistic
interpretation of harmonic functions as solutions of Dirichlet
problems is not true anymore. The following result, which
follows easily from Theorem \ref{2l}, provides some probabilistic
interpretation to these kind of harmonic functions.

\begin{thm}\label{2m}
Suppose that $D$ is a bounded Lipschitz domain and $\mu$ is
a finite measure on $\partial D$. Define
\[
h(x)=\int_{\partial D}M_D(x, z)\mu(dz).
\]
Then for any Borel measurable subset $A\subset \partial D$,
\[
P_x^h(\lim_{t\uparrow\tau_D}X(t)\in A)=\frac1{h(x)}\int_A
M_D(x, z)\mu(dz).
\]
\end{thm}

In particular, when $D=B(O, r)$ and
$h(x)=\int_{\partial D}M_D(x, z)\sigma(dz)$, where $\sigma$ is the
surface measure, $\lim_{t\uparrow\tau_D}X(t)$ is distributed uniformly
on $\partial D$ under $P_0^h$.

\vspace{.1in}

\section{Integral Representations of Positive Harmonic Functions}

Functions which are (super)harmonic in $D$ with respect to $X$ are studied
in Landkof \cite{ns:lan} and  Bogdan \cite{kr:bog}. 
However, it seems that no one has studied the integral
representations of this kind of (super)harmonic functions. We intend
to establish such a representation. To prove the uniqueness
of such a representation
theorem we need the following result:

\begin{lemma}\label{3a}
Suppose that $D$ is a bounded Lipschitz domain. If a function $f$
satisfies the following
\begin{equation}
\int_{D^c}\frac{f(z)}{|y-z|^{n+\alpha}}dz=0, \qquad \forall y\in D,
\end{equation}
Then $f=0$ almost everywhere on $D^c$.
\end{lemma}

\pf Without loss of generality we can assume that the origin $O$
is in $D$.

We claim that for all $0\le m\le k$,
\begin{equation}
\int_{D^c}\prod_{j=1}^m(x_{i_j}-y_{i_j})|x-y|^{-(n+\alpha+2k)}f(y)dy=0,
\qquad x\in D,
\end{equation}
where for each $j$, $1\le i_j\le n$. We are going to prove the claim 
by induction on $k$. 

The case of $k=0, m=0$ follows by the assumption. Take the partial
derivative of (4.1) with respect to $x_j$ we get
\[
0=\frac{\partial}{\partial x_j}\int_{D^c}|x-y|^{-n-\alpha}f(y)dy
=-(n+\alpha)\int_{D^c}(x_j-y_j)|x-y|^{-n-\alpha-2}f(y)dy,
\]
Thus
\begin{equation}
\int_{D^c}(x_j-y_j)|x-y|^{-n-\alpha-2}f(y)dy=0, \qquad x\in D,
\end{equation}
that is, the claim is true for $k=m=1$. Now take the partial
derivative of (4.3) with respect to $x_j$ we get
\begin{eqnarray*}
0&=&\frac{\partial}{\partial x_j}
\int_{D^c}(x_j-y_j)|x-y|^{-n-\alpha-2}f(y)dy\\
&=&\int_{D^c}|x-y|^{-n-\alpha-2}f(y)dy
-(n+\alpha+2)\int_{D^c}(x_j-y_j)^2|x-y|^{-n-\alpha-4}f(y)dy
\end{eqnarray*}
Summing the above from $j=1$ to $j=n$ we get
\[
(\alpha+2)\int_{D^c}|x-y|^{-n-\alpha-2}f(y)dy=0,
\]
which implies that the claim is true for the case of $k=1, m=0$.
Therefore the claim is true for $k=1$.

Now we assume that the claim is true for all $0\le m\le k\le N$. 
Take the partial derivative of (4.2) with respect to $x_{i_{m+1}}$ we get
\begin{eqnarray*}
0&=&\frac{\partial}{\partial x_{i_{m+1}}}
\int_{D^c}\prod_{j=1}^m(x_{i_j}-y_{i_j})|x-y|^{-(n+\alpha+2k)}f(y)dy\\
&=&\int_{D^c}\frac{\partial}{\partial x_{i_{m+1}}}\left(
\prod_{j=1}^m(x_{i_j}-y_{i_j})|x-y|^{-(n+\alpha+2k)}\right)f(y)dy\\
&=&-(n+\alpha+2k)\int_{D^c}\prod_{j=1}^{m+1}
(x_{i_j}-y_{i_j})
|x-y|^{-(n+\alpha+2(k+1))}f(y)dy\\
&&+\int_{D^c}\frac{\partial}{\partial x_{i_{m+1}}}\left(
\prod_{j=1}^m(x_{i_j}-y_{i_j})\right)|x-y|^{-(n+\alpha+2k)}f(y)dy\\
&=&-(n+\alpha+2k)\int_{D^c}\prod_{j=1}^{m+1}
(x_{i_j}-y_{i_j})
|x-y|^{-(n+\alpha+2(k+1))}f(y)dy,
\end{eqnarray*}
where in the last equality we used the induction assumption. Therefore
the claim is true for all $0\le m\le k\le N+1$, and hence the claim 
is always true. 

Evaluate (4.2) at $x=O$ we get that for any non-negative integer $k$,
and any multi--index $\beta=(\beta_1, \cdots, \beta_n)$ with $|\beta|
=\beta_1+\cdots+\beta_n\le k$,
\[ 
\int_{D^c}y^{\beta}|y|^{-2k}\frac{f(y)}{|y|^{n+\alpha}}dy=0
\]
where $y^{\beta}=y_1^{\beta_1}\cdots y_n^{\beta_n}$. Since the linear
span of the set $\{y^{\beta}|y|^{-2k}: |\beta|\le k\}$ is an algebra
of real--valued continuous functions on $D^c$ which separates
points in $D^c$ and vanishes at infinity,
by the Stone--Weierstrass Theorem the linear span of
$\{y^{\beta}|y|^{-2k}: |\beta|\le k\}$ is dense in 
$C_{\infty}(D^c)$ with respect to the uniform topology.
Here $C_{\infty}(D^c)$ is the space of continuous functions on $D^c$ which vanishes at infinity. Thus for all $\phi\in
C_{\infty}(D^c)$,
\[
\int_{D^c}\phi(y)\frac{f(y)}{|y|^{n+\alpha}}dy=0
\]
which implies that $f(y)|y|^{-(n+\alpha)}=0$ almost everywhere on $D^c$.
Therefore $f=0$ almost everywhere on $D^c$.
\qed

\vspace{.1in}

\begin{thm}\label{3d}
Suppose that $D$ is a bounded
domain in $\RR^n$. If $h$ and $f$ are both harmonic
in $D$ with respect to $X$ with $h=f$ in $D$, then $h=f$ in $\RR^n$.
\end{thm}

\pf
Take $x_0\in D$ and $B(x_0, r)\subset\overline{B(x_0, r)}\subset D$, then 
it follows from Theorem \ref{1b} that for any $x\in B(x_0, r)$,
\[
E_x[h(X_{\tau_{B(x_0, r)}})]=h(x)=f(x)=E_x[f(X_{\tau_{B(x_0, r)}})].
\]
Therefore we have
\[
E_x[(h-f)(X_{\tau_{B(x_0, r)}})]=0, \qquad x\in B(x_0, r).
\]
By Theorem 1.4 of Chen and Song \cite{s1:che} we know that
for all $x\in B(x_0, r)$,
\begin{eqnarray*}
E_x[(h-f)(X_{\tau_{B(x_0, r)}})]&=&
\int_{B(x_0, r)^c}K_{B(x_0, r)}(x, z)(h-f)(z)dz\\
&=& A(n, \alpha)\int_{B(x_0, r)^c}\left(\int_{B(x_0, r)}
G_{B(x_0, r)}(x, y)\frac1{|y-z|^{n+\alpha}}dy\right)(h-f)(z)dz\\
&=&  A(n, \alpha) 
\int_{B(x_0, r)}G_{B(x_0, r)}(x, y)\left(\int_{B(x_0, r)^c}
\frac{(h-f)(z)}{|y-z|^{n+\alpha}}dz\right)dy.
\end{eqnarray*}
Therefore by general potential theory (see Section 5.2 of 
\cite{kl:chu}, for instance) we know that the function
\[
y\mapsto\int_{B(x_0, r)^c}\frac{(h-f)(z)}{|y-z|^{n+\alpha}}dz
\]
is zero almost everywhere on $B(x_0, r)$.
Since the function above
is continuous in $B(x_0, r)$, we have
\[
\int_{B(x_0, r)^c}\frac{(h-f)(z)}{|x-z|^{n+\alpha}}dz=0, \qquad \forall x\in 
B(x_0, r).
\]
It follows from Lemma \ref{3a} we know that $u-f=0$
almost everywhere on $B(x_0, r)^c$, and the proof is finished..
\qed

\vspace{.1in}

\noindent{\bf Remark 4.1}
In fact the proof actually shows that for
a function $h$ 
harmonic in a domain $D$, the values of $h$ in
any ball $B(x_0, r)\subset\overline{B(x_0, r)}\subset D$ determine
$h$ uniquely.

\begin{thm}\label{3b}
Suppose $D$ is a bounded Lipschitz domain in $\RR^n$.
If $f$ is a non-negative harmonic function in $D$ 
with respect to $X$, then
there exists a unique finite measure $\mu$ on $\partial D$ such that
the restriction of $f$ to $D$ can be written as
\begin{equation}\label{eqn:4.4}
f(x)=\int_{D^c}K_D(x, z)f(z)dz+\int_{\partial D}M_D(x, z)\mu(dz),
\qquad \forall x\in D.
\end{equation}
\end{thm}

\pf First we are going to show that the difference
\[
f(x)-\int_{D^c}K_D(x, z)f(z)dz
\]
is non-negative. 
Take a sequence of domains $D_m$ such that $D_m\subset
\overline{D_m}\subset D_{m+1}\subset\overline{D_{m+1}}
\subset D$ and set $\tau_m=\tau_{D_m}$. Then $\tau_m\uparrow
\tau_D$. For any $x\in D$, since $P_x(X_{\tau_{D-}}\neq X_{\tau_D})=1$,
we know that $P_x(\tau_D=\tau_m\, \mbox{ for some $m\ge 1$})
=1$.
From Theorem \ref{1b} we know that
\begin{eqnarray*}
f(x)&=& E_xf(X_{\tau_m})\\
&=& E_x[f(X_{\tau_D}); \tau_m=\tau_D]+E_x[f(X_{\tau_m}); \tau_m<\tau_D]\\
&\ge& E_x[f(X_{\tau_D}); \tau_m=\tau_D].
\end{eqnarray*}
Therefore
\[
f(x)\ge E_x[f(X_{\tau_D})].
\]
Therefore by Theorem \ref{2d} we know that there exists a unique
finite measure
$\mu$ on $\partial D$ such that
\[
f(x)-E_x[f(X_{\tau_D})]=\int_{\partial D}M_D(x, z)\mu(dz),
\qquad \forall x\in D.
\]
\qed

\vspace{.1in}
 
From the above theorem we can easily get the following

\begin{thm}\label{3c}
If $D$ is a bounded Lipschitz domain, then the restriction to
$D$ of any non-negative function $f$ which is superharmonic in $D$
with respect to $X$ can be written as
\[
f(x)=\int_{D^c}K_D(x, z)f(z)dz+\int_DG_D(x, y)\nu(dy)+
\int_{\partial D}M_D(x, z)\mu(dz),
\]
where $\nu$ and $\mu$ are finite measures on $D$ and
$\partial D$ respectively.
\end{thm}

\pf Similar to the first part of the proof of Theorem \ref{3b},
we have that the function
\[
f(x)-E_x[f(X_{\tau_D})]
\]
is a non-negative function which 
vanishes outside $D$ and is superharmonic 
in $D$ with respect to $X$. Hence it is a non-negative 
function which is 
harmonic in $D$ with respect to $X^D$. Now our claim follows from Theorem 
\ref{2f}.
\qed

\vspace{.1in}

However, the above decomposition is not unique
anymore. This non-uniqueness is due to the
following relation between the Green function $G_D$
and Poisson kernel $K_D$ established in  Theorem 1.4 of Chen and Song
\cite{s1:che}:
\[
K_D(x, z)=A(n, \alpha)\int_D\frac{G_D(x, y)}{|y-z|^{n+\alpha}}dy.
\]
Using this relation, we can 
absorb the first term on the right hand side above, or part of it,
into the second term.

\vspace{.1in}

\section{Extensions}

The results of this paper can be extended to dimension $n=1$,
by noting that Green functions and Poisson kernels on
bounded open intervals for
symmetric stable processes are explicitly known (see, for example,
\cite{g1:blu}). In particular
(\ref{eqn:2.1})--(\ref{eqn:2.2}) hold
for bounded open interval $U$ in $\RR$. 
Note that the Kato class $\KK_{1, \alpha}$ is defined by (\ref{eqn:kato})
for $0<\alpha <2$
(cf. Corollary 4 of \cite{g1:blu} and \cite{zhao}).

\vspace{.5in}
\begin{singlespace}

\end{singlespace}

\end{doublespace}

\end{document}